\documentstyle[12pt]{article}
\input amssym.def
\input amssym

\def\sign{\mathop{\rm sign}\nolimits}
\begin{document}
\title{Zeros of the hypergeometric polynomial $F(-n,b;\,c;\,z)$}
\author{K. Driver\thanks{Research of the first author is supported
by the John Knopfmacher Centre for Applicable Analysis and
Number Theory, University of the Witwatersrand.} \ and K. Jordaan}
\date{}
\maketitle
\smallskip

\begin{center}
School of Mathematics\\  University of the
Witwatersrand\\  Johannesburg, South Africa
\end{center}
\medskip

\begin{abstract}
Our interest lies in describing the zero behaviour of Gauss
hypergeometric polynomials $F(-n,b;\,c;\,z)$ where $b$ and $c$
are arbitrary parameters. In general, this problem has
not been solved and even when $b$ and $c$ are both real, the only
cases that have been fully analyzed impose additional restrictions
on $b$ and $c$. We review recent results that have been proved
for the zeros of several classes of hypergeometric polynomials
$F(-n,b;\,c;\,z)$ where $b$ and $c$ are real. We show that
the number of real zeros of $F(-n,b;\,c;\,z)$ for arbitrary real values
of the parameters $b$ and $c$, as well as the intervals in which
these zeros (if any) lie, can be deduced from corresponding results
for Jacobi polynomials.
\end{abstract}
\bigskip

\noindent AMS MOS Classification:\quad  33C05, 30C15.
\medskip


\section{Introduction}
The Gauss hypergeometric function, or ${}_2F_1$, is defined by
\[
F(a,b;\, c;\, z) = 1+\sum_{k=1}^\infty \frac{(a)_k(b)_k}{(c)_k}\,
\frac{z^k}{k!}, \quad |z| < 1,
\]
where $a,\,b$ and $c$ are complex parameters and
\[
(\alpha)_k =\alpha(\alpha+1)\dots(\alpha+k-1) = \Gamma(\alpha+k)
\big/ \Gamma(\alpha)
\]
is Pochhammer's symbol. When $a=-n$ is a negative integer, the series
terminates and reduces to a polynomial of degree $n$, called a
hypergeometric polynomial. Our focus lies in the location of the zeros
$F(-n,b;\,c;\,z)$ for real values of $b$ and $c$.

Hypergeometric polynomials are connected with several different types
of orthogonal polynomials, notably Chebyshev, Legendre, Gegenbauer
and Jacobi polynomials. In the cases of Chebyshev and Legendre
polynomials, the connection demands fixed special values of the
parameters $b$ and $c$, namely, (cf. [1], p.561)
\[
F\left(-n,n;\, \frac 12;\,z\right) = T_n(1-2z)
\]
and
\[
F\left(-n,n+1;\, 1;\,z\right) = P_n(1-2z),
\]
respectively. However, in the cases of Gegenbauer and Jacobi polynomials,
we have
\begin{equation}
F\left(-n,n+2\lambda;\, \lambda+\frac 12;\,z\right) = \frac {n!}
{(2\lambda)_n} C_n^\lambda(1-2z)
\end{equation}
and
\begin{equation}
F\left(-n,\alpha+\beta+1+n;\, \alpha+1;\,z\right) = \frac {n!}
{(\alpha+1)_n} {\cal P}_n^{(\alpha,\beta)}(1-2z),
\end{equation}
respectively. Since the zeros of orthogonal polynomials are well
understood, we expect the connections (1.1) and (1.2) to be very
useful in analyzing the zeros of $F(-n,b;\, c;\, z)$. Conversely,
if the zeros of $F(-n,b;\, c;\, z)$ are known, this leads to new
information about the zero distribution of Gegenbauer or Jacobi
polynomials for values of their parameters that lie outside the range
of orthogonality of these polynomials.

This paper is organized as follows. In Section 2 we give a self-contained
review of recent results regarding the zeros of several special classes
of hypergeometric polynomials. Section 3 contains results
originally due to Klein [9] which
detail the numbers and location of real zeros of $F(-n,b;\,c;\, z)$
for arbitrary real values of $b$ and $c$.
We provide simple proofs using results proved in [13].

\section{Zeros of special classes of hypergeometric polynomials}
\setcounter{equation}{0}

We begin with a few general remarks. Since we shall assume throughout
our discussion that $b$ and $c$ are real parameters, we know that
all zeros of
$F(-n,b;\,c;\, z)$ must occur in complex conjugate pairs. In particular,
if $n$ is odd, $F$ must always have at least one real zero. Further,
if $b=-m$ where $m<n$, $m\in\Bbb N$, $F(-n,b;\,c;\, z)$ reduces to a
polynomial of degree $m$. However, since we are interested in the
behaviour of the zeros of $F(-n,b;\,c;\, z)$ as $b$ and/or $c$
vary through real values, we shall adopt the convention that
$F(-n,-m;\,c;\, z)= \lim_{b\to -m}F(-n,b;\,c;\, z)$. This ensures
that the zeros of $F$ vary continuously with $b$ and $c$. Note also
that $F(-n,b;\,c;\, z)$ is not defined when $c=0,-1,\dots,-n+1$.
Regarding the multiplicity of zeros, a hypergeometric function
$w=F(a,b;\,c;\,z)$ satisfies the differential equation
\[
z(1-z)w''+\big[c-(a+b+1)z\big] w'-abw =0,
\]
so if $w(z_0)=w'(z_0)=0$ at some point $z_0\ne 0$ or 1,
it would follow that $w\equiv 0$. Thus multiple zeros of $F(-n,b;\,c;\,z)$
can only occur at $z=0$ or 1.

\subsection{Quadratic transformations}

The class of hypergeometric polynomials that admit a quadratic transformation
is specified by a necessary and sufficient condition due to Kummer
(cf.\ [1], p.560). There are twelve polynomials in this class
(cf.\ [14], p.124)
\[
\begin{array}{lll}
F \left( -n, b;\, 2b;\, z \right)
& F\left( -n,b;\, -n-b+1;\, z\right)
& F\left( -n,b;\, \frac{-n+b+1}{2};\, z\right)
\\[.15in]
F \left( -n, b;\, \frac 12;\, z \right)
& F\left( -n,-n+\frac 12;\, c;\, z\right)
& F\left( -n,b;\, -n+b+\frac 12;\, z\right)
\\[.15in]
F \left( -n, b;\, \frac 32;\, z \right)
& F\left( -n,-n-\frac 12;\, c;\, z\right)
& F\left( -n,b;\, -n+b-\frac 12;\, z\right)
\\[.15in]
F \left( -n, b;\, -2n;\, z \right)
&F\left( -n,b;\, b+n+1;\, z\right)
&F\left( -n,n+1;\, c;\, z\right)\,.
\end{array}
\]
The most important polynomial in this class is $F(-n,b;\,2b;\,z)$
because complete analysis of its zero distribution for all real values
of $b$ (cf.\ [4], [5]) leads to corresponding results for the zeros of the
Gegenbauer polynomials $C_n^\lambda(z)$ for all real values of the
parameter $\lambda$ (cf.\ [6]).
\bigskip

\noindent{\bf Theorem 2.1. } {\it Let $F=F(-n,b;\, 2b;\, z)$ where
$b$ is real.
\begin{itemize}
\item[(i)] For $b> -\frac 12$, all zeros of $F(-n,b;\,2b;\,z)$ are
simple and lie on the circle $|z-1|=1$.

\item[(ii)] For $-\frac 12 -j < b< \frac 12 -j$, $j=1,2,\dots
\left[ \frac n2\right]-1$, $(n-2j)$ zeros of $F$ lie on the circle
$|z-1|=1$. If $j=2k$ is even, there are $k$ non-real zeros of $F$
in each of the four regions bounded by the circle $|z-1|=1$ and
the real axis. If $j=2k+1$ is odd, there are $k$ non-real zeros of
$F$ in each of the four regions described above and the remaining two
zeros are real.

\item[(iii)] If $n$ is even, for $-\left[\frac n2\right] <b<-\left[
\frac n2\right] + \frac 12$, no  zeros of $F$ lie on $|z-1|=1$.
If $n=4k$, all zeros of $F$ are non-real whereas if $n=4k+2$, two
zeros of $F$ are real and $4k$ are non-real. If $n$ is odd, for
$-1-\left[\frac n2\right] <b<-\left[\frac n2\right] +\frac 12$, only
the fixed real zero of $F$ at $z=2$ lies on $|z-1|=1$. If $n=4k+1$,
$n-1=4k$ zeros of $F$ are non--real whereas if $n=4k+3$, two further
zeros are real and the remaining $4k$ are non--real.

\item[(iv)] For $j-n<b<j-n+1$, $j=1,2,\dots\left[\frac n2\right]-1$,
$(n-2j)$ zeros of $F$ are real and greater than 1. If $j=2k$ is even,
all remaining $2j$ zeros of $F$ are non--real with $k$ zeros in
each of the regions described above; while if $j=2k+1$, $4k$ zeros
are non--real as before and 2 are real.

\item[(v)] For $b<1-n$, all zeros of $F(-n,b;\,2b;\,z)$ are real and
greater than 1. As $b\to-\infty$, all the zeros of $F$ converge to
the point $z=2$.
\end{itemize} }
\medskip

An analogous theorem which describes the behaviour of the zeros of
$C_n^\lambda(z)$ can be found in [6], Section 3 or [7], Theorem 1.2.

For the polynomial $F\left(-n,b;\,\frac 12;\, z\right)$
the following result has been proved in [7], Theorem 2.3.
\bigskip

\noindent{\bf Theorem 2.2. } {\it Let $F=F\left(-n,b;\,\frac 12;\, z\right)$
with $b$ real.
\begin{itemize}
\item[(i)] For $b>n-\frac 12$, all $n$ zeros of $F$ are real and simple
and lie in $(0,1)$.
\item[(ii)] For $n-\frac 12-j<b<n+\frac 12 -j$, $j=1,2,\dots,n-1$,
$(n-j)$ zeros of $F$ lie in $(0,1)$ and the remaining $j$ zeros
of $F$ form $\left[\frac j2\right]$ non-real complex pairs of zeros
and one real zero lying in $(1,\infty)$ when $j$ is odd.
\item[(iii)] For $0<b<\frac 12$, $F$ has $\left[\frac n2\right]$ non-real
complex conjugate pairs of zeros with one real zero in $(1,\infty)$
when $n$ is odd.

\item[(iv)] For $-j<b<-j+1$, $j=1,2,\dots,n-1$, $F$ has exactly $j$
real negative zeros. There is exactly one further real zero greater
than 1 only when $(n-j)$ is odd and all the remaining zeros of $F$
are non--real.

\item[(v)] For $b<1-n$, all zeros of $F$ are real and negative and
converge to zero as $b\to-\infty$.
\end{itemize} }
\bigskip

A very similar theorem is proved for the zeros of $F\left( -n,b;\,
\frac 32;\, z\right)$ in [7], Theorem 2.4 with only minor differences
of detail.

For the hypergeometric polynomial $F(-n,b;\, -2n;\, z)$, less complete
results have been proved. We have (cf.\ [8] Theorem 3.1 and
Corollary 3.2) the following.
\bigskip

\noindent{\bf Theorem 2.3. } {\it Let $F=F\left(-n,b;\, -2n;\, z\right)$
with $b$ real.
\begin{itemize}
\item[(i)] For $b>0$, $F$ has $n$ non-real zeros if $n$ is even whereas
if $n$ is odd, $F$ has exactly one real negative zero and the remaining
$(n-1)$ zeros of $F$ are all non-real.
\item[(ii)] For $-n<b<0$, if $-k<b<-k+1$, $k=1,\dots,n$,
$F$ has $k$ real zeros in the interval $(1,\infty)$. In addition,
if $(n-k)$ is even, $F$ has $(n-k)$ non-real zeros whereas if $(n-k)$
is odd, $F$ has one real negative zero and $(n-k-1)$ non-real zeros.

\item[(iii)] For $-n>b>-2n$, if $-n-k>b>-n-k-1$, $k=0,1,\dots,n-1$,
$F$ has $(n-k)$ real zeros in the interval $(1,\infty)$.
In addition, if $k$ is even $F$ has $k$ non-real zeros while if $k$
is odd, $F$ has one real zero in $(0,1)$ and $(k-1)$ non-real zeros.

\item[(iv)] For $b<-2n$, all $n$ zeros of $F$ are non--real for $n$
even whereas for $n$ odd, $F$ has exactly one real zero in the
interval $(0,1)$.
\end{itemize} }
\medskip

The identities (cf.\ [7], Lemma 2.1)
\begin{equation}
F(-n,b;\, c;\, 1-z) = \frac{(c-b)_n}{(c)_n} F(-n,b;\, 1-n+b-c;\,z)
\end{equation} 
and
\begin{equation}
F(-n,b;\, c;\, z) = \frac{(b)_n}{(c)_n}(-z)^n F\left(-n,1-c-n;\,
1-b-n;\,\frac 1z\right)
\end{equation} 
hold for $b$ and $c$ real, $c\ne \{0,-1,\dots,-n+1\}$.
Applying (2.1) and (2.2) to each of the polynomials $F(-n,b;\,2b;\,z)$,
$F\left(-n,b;\,\frac 12;\, z\right)$, $F\left(-n,b;\,\frac 32;\,z
\right)$ and $F(-n,b;\,-2n;\,z)$ in turn, we obtain the remaining
eight polynomials in the quadratic class. It is then an easy task
to deduce analogous results for their zero distribution.

A similar set of results has been proved for the sixteen hypergeometric
polynomials in the cubic class. Again, this class arises from a necessary
and sufficient condition (cf.\ [2], p.67) and details can be found
in [7].

\section{The real zeros of $F(-n,b;\, c;\,z)$ for $b$ and $c$ real} 

The results proved below are due to Klein [9] who considered the zeros
of more general hypergeometric functions (not necessarily polynomials).
Klein's proof is geometric and difficult to penetrate. A more transparent
perspective in the polynomial case may be provided by the approach
given here.

The classical equation linking the hypergeometric polynomial
$F(-n,b;\,c;\, z)$  with Jacobi polynomials ${\cal P}_n^{(\alpha,\beta)}(z)$
is given by (1.2). We will find an alternative expression
(cf.\ [12], p.464, eqn. (142))
\setcounter{equation}{0}
\begin{equation}
F(-n,b;\, c;\,z) = \frac{n!z^n}{(c)_n} {\cal P}_n^{(\alpha,\beta)}
\left( 1-\frac 2z \right),
\end{equation}    
where $\alpha =-n-b$ and $\beta = b-c-n$, more suited to our analysis.
The number of real zeros
of ${\cal P}_n^{(\alpha,\beta)}(x)$ in the intervals $(-1,1)$, $(-\infty,1)$
and $(1,\infty)$ are given by the Hilbert-Klein formulas (cf.\ [13],
p.145, Theorem 6.72), also known to Stieltjes.
We use Klein's symbol
\[
E(u) = \left\{ \begin{array}{ll}
0 & \mbox{if } u\le 0\\ {[u]} & \mbox{if $u>0$, $u \ne$ integer}\\
u-1 & \mbox{if } u=1,2,3,\dots
\end{array} \right.\quad .
\]
\bigskip

Noting that under the linear fractional transformation $w=1-\frac 2z$,
the intervals $1<w<\infty$, $-\infty<w<-1$ and $-1<w<1$ correspond to
$-\infty<z<0$, $0<z<1$ and $1<z<\infty$ respectively, we can use
equation (3.1) to rephrase the Hilbert-Klein formulas for hypergeometric
polynomials.
\bigskip

\noindent{\bf Theorem 3.1. } {\it Let $b,\, c\in\Bbb R$ with $b,c,
c-b\ne 0$, $-1,\dots,-n+1$. Let
\begin{eqnarray}
X &=& E\left\{ \frac 12 \left( |1-c| - |n+b| - |b-c-n| +1\right) \right\}
\\ 
Y &=& E\left\{ \frac 12 \left( -|1-c| + |n+b| - |b-c-n| +1\right) \right\}
\\ 
Z &=& E\left\{ \frac 12 \left( -|1-c| - |n+b| + |b-c-n| +1\right) \right\}.
\end{eqnarray}
Then the numbers of zeros of $F(-n,b;\,c;\, z)$ in the intervals
$(1,\infty)$, $(0,1)$ and $(-\infty,0)$ respectively are
\begin{eqnarray}
N_1 = \left\{ \begin{array}{ll}
2[(X+1)/2] &\mbox{if } (-1)^n {-b \choose n} {b-c \choose n} >0\\[.10in]
2[X/2]+1 &\mbox{if } (-1)^n {-b \choose n} {b-c \choose n}<0
\end{array} \right.
\\[.15in] 
N_2 = \left\{ \begin{array}{ll}
2[(Y+1)/2] &\mbox{if } {-c \choose n} {b-c \choose n} >0\\[.10in]
2[Y/2]+1 &\mbox{if } {-c \choose n} {b-c \choose n}<0
\end{array} \right.
\\[.15in]     
N_3 = \left\{ \begin{array}{ll}
2[(Z+1)/2] &\mbox{if } {-c \choose n} {-b \choose n} >0\\[.10in]
2[Z/2]+1 &\mbox{if } {-c \choose n} {-b \choose n}<0.
\end{array} \right.
\end{eqnarray} } 
\bigskip

\noindent{\bf Proof. } The expressions all follow immediately from the
Hilbert-Klein formulas (cf.\ [13], p.145, Thm.\ 6.72)
together with equation (3.1).
$\Box$
\bigskip

\noindent{\bf Theorem 3.2. } {\it Let $F=F(-n,b;\, c;\, z)$ where
$b,\,c\in\Bbb R$ and $c>0$.
\begin{itemize}
\item[(i)] For $b>c+n$, all zeros of $F$ are real and lie in the interval
$(0,1)$.
\item[(ii)] For $c<b<c+n$, $c+j-1<b<c+j$, $j=1,2,\dots,n$; $F$ has $j$
real zeros in $(0,1)$. The remaining $(n-j)$ zeros of $F$ are all non-real
if $(n-j)$ is even while if $(n-j)$ is odd, $F$ has $(n-j-1)$ non-real
zeros and one additional real zero in $(1,\infty)$.
\item[(iii)] For $0<b<c$, all the zeros of $F$ are non-real if $n$
is even, while if $n$ is odd, $F$ has one real zero in $(1,\infty)$
and the other $(n-1)$ zeros are non-real.
\item[(iv)] For $-n<b<0$, $-j<b<-j+1$, $j=1,2,\dots,n$, $F$ has $j$
real negative zeros.  The remaining $(n-j)$ zeros of $F$ are all
non-real if $(n-j)$ is even, while if $(n-j)$ is odd, $F$ has $(n-j-1)$
non-real zeros and one additional real zero in $(1,\infty)$.
\item[(v)] For $b<-n$, all zeros of $F$ are real and negative.
\end{itemize} }
\bigskip

\noindent{\bf Proof. } We use the identity (cf.\ [1], p.559, (15.3.4))
\begin{equation}
F(-n,b;\, c;\, z) = (1-z)^n F\left(-n,c-b;\, c;\, \frac z{z-1} \right)
\end{equation} 
to show that (i) $\Rightarrow$ (v) and (ii) $\Rightarrow$ (iv) so that
it will suffice to prove (i), (ii) and (iii) above.
\newline (i) $\Rightarrow$ (v): If $b<-n$ then $c-b>c+n$ and
by (i), all zeros of $F(-n,c-b;\, c;\, w)$ are real and lie in the
interval $(0,1)$. Since $w=\frac z{z-1}$ maps $(-\infty, 0)$ to
$(0,1)$, (v) follows from (3.8).\newline
(ii) $\Rightarrow$ (iv): If $-j<b<-j+1$, $j=1,2,\dots,n$, then
$c+j-1<c-b<c+j$, $j=1,2,\dots,n$. By (ii), since $w=\frac z{z-1}$
maps $(-\infty,0)$ to $(0,1)$ and $(1,\infty)$ to $(1,\infty)$,
(iv) follows again from (3.8). To prove (i), (ii) and (iii), we note
that in each part, $b>0$ (and of course $c>0$ by assumption). Then
\begin{equation}
\sign{-b\choose n} =(-1)^n, \quad \sign
{-c\choose n} =(-1)^n.
\end{equation} 
\begin{itemize}
\item[(i)] Suppose $b>c+n$. Then $b-c>n$ and
\begin{equation}
\sign {b-c\choose n} >0 \mbox{ for all } n.
\end{equation} 
Considering (3.5), (3.6) and (3.7) with (3.9) and (3.10), we observe
that
\begin{eqnarray*}
N_1 &=& 2\left[(X+1)/2\right], \quad N_3 = 2\left[(Z+1)/2\right],
\\[.15in]
N_2 &=& \left\{ \begin{array}{ll}
2\left[(Y+1)/2\right] &\mbox{for $n$ even}\\[.10in]
2\left[Y/2\right] + 1 &\mbox{for $n$ odd} \end{array} \right. .
\end{eqnarray*}
Assume now that $c>1$. Then for $b>c+n$, we have from (3.2), (3.3)
and (3.4) that $X=0$, $Y=n$, $Z=0$. Substituting these values into
$N_1$, $N_2$ and $N_3$ yields the result. A similar calculation shows
that the same result is obtained when $0<c<1$.

\item[(ii)] For $c+j-1<b<c+j$, $j=1,2,\dots,n$, we find that
$\sign{b-c\choose n}=(-1)^{n-j}$. Then from (3.5), (3.6),
(3.7) we see that
\begin{eqnarray*}
N_1 &=& \left\{ \begin{array}{ll}
2\left[(X+1)/2\right] &\mbox{for $(n-j)$ even}\\[.10in]
2\left[X/2\right] + 1 &\mbox{for $(n-j)$ odd} \end{array} \right. ,
\\[.15in]
N_2 &=& \left\{ \begin{array}{ll}
2\left[(Y+1)/2\right] &\mbox{for $j$ even}\\[.10in]
2\left[Y/2\right] + 1 &\mbox{for $j$ odd} \end{array} \right. ,
\\[.15in]
N_3 &=& 2\left[(Z+1)/2\right].
\end{eqnarray*}
It follows from (3.2), (3.3) and (3.4) by an easy calculation that
$X=0$, $Y=j$, $Z=0$ and we deduce that $N_1=\left\{ \begin{array}{ll}
0 &\mbox{if $(n-j)$ is even}\\ 1 &\mbox{if $(n-j)$ is odd} \end{array}
\right.$, $N_2=j$ and $N_3=0$ which proves (ii).

\item[(iii)] For $0<b<c$, $\sign {b-c\choose n} = (-1)^n$.
Then $N_1 =\left\{\begin{array}{ll} 2\left[(X+1)/2\right] &\mbox{if
$n$ is even}\\2\left[X/2\right]+1 &\mbox{if $n$ is odd} \end{array}
\right.,$ $N_2 =2\left[(Y+1)/2\right]$, $N_3 = 2\left[(Z+1)/2\right]$.
Also, we find $X=0$, $Y=0$ and $Z=0$ which completes the proof of (iii)
and hence the theorem.
\end{itemize}

For $c<0$, the range of values of $b$ and $c$ that have to be considered
can be reduced if we use the identities (2.1) and (2.2). Since the
real zeros of $F(-n,b;\, c;\, z)$ are now known for all $c>0$ and
$b\in\Bbb R$ from Theorem 3.2, it follows from (2.1) that we need only
consider $c-b>1-n$. Similarly, from (2.2) and Theorem 3.2, we can assume
$b>1-n$. We split the result for $c<0$ into the cases where $b>0$
and $1-n<b<0$.
\bigskip

\noindent{\bf Theorem 3.3. } {\it Let $F=F(-n,b;\, c;\, z)$. Suppose
that $c<0$, $b>0$, $c-b>1-n$. Then
\begin{itemize}
\item[(i)] $1-n<c-b<0$ and $0<b<n-1$ and $1-n<c<0$.
\item[(ii)] If $-k<c<-k+1$, $k=1,\dots,n-1$ and
\[
-j<c-b<-j+1, \quad j=1,\dots,n-1,
\]
then $F(-n,b;\,c;\,z)$ has $(j-k)\ge 0$ real zeros in $(0,1)$. For the
remaining $(n-j+k)$ zeros of $F$
\begin{itemize}
\item[(a)] $(n-j+k)$ are non-real if $(n-j)$ and $k$ are even
\item[(b)] $(n-j+k-1)$ are non-real and one real zero lies in $(1,\infty)$
if $(n-j)$ is odd and $k$ is even
\item[(c)] $(n-j+k-1)$ are non-real if $(n-j)$ is even, $k$ odd and one
zero is real and negative
\item[(d)] $(n-j+k-2)$ are non-real if $(n-j)$ is odd and $k$ is odd
with one real negative zero and one real zero in $(1,\infty)$.
\end{itemize}
\end{itemize} }
\bigskip

\noindent{\bf Proof. }
\begin{itemize}
\item[(i)] This follows immediately from $c<0$, $b>0$, $c-b>1-n$.
\item[(ii)] For $c<0$, $b>0$, $c-b>1-n$, we have
\[
|1-c|=1-c, \quad |b+n|=b+n, \quad |b-c-n|=c-b+n
\]
and it follows from (3.2), (3.3) and (3.4) that
\[
X=E(1-c-n), \quad Y=E(b), \quad Z=E(c-b).
\]
Since $1-c-n<0$ and $c-b<0$, $X=Z=0$. Now $\sign {-b\choose n}
=(-1)^n$ and for $k=1,\dots,n-1$, $-k<c<-k+1 \Rightarrow \sign
{-c\choose n} =(-1)^{n-k}$, while for $-j<c-b<-j+1$, $j=1,
\dots,n-1$, $\sign {b-c\choose n} =(-1)^{n-j}$.
Therefore, from (3.5), (3.6) and (3.7),
\begin{eqnarray}
N_1 &=& \left\{ \begin{array}{ll}
0 &\mbox{if $(n-j)$ even}\\ 1 &\mbox{if $(n-j)$ odd}
\end{array}\right.
\\ 
N_2 &=& \left\{ \begin{array}{ll}
2\left[(Y+1)/2\right] &\mbox{if $(j-k)$ is even}\\
2\left[Y/2\right]+1 &\mbox{if $(j-k)$ is odd} \end{array}\right. ,
\quad Y=E(b)
\\ 
N_3 &=& \left\{ \begin{array}{ll}
0 &\mbox{if $k$ even}\\ 1 &\mbox{if $k$ odd} \end{array}\right. .
\end{eqnarray} 
Now for $j>b-c>j-1$ and $-k<c<-k+1$, $b\in (j-k-1,\,j-k+1)$,
$j-k=1,2,\dots,n-2$. If $b\in (j-k-1,\,j-k)$, $Y=E(b) =j-k-1$,
whereas if $b\in (j-k,\,j-k+1)$, $Y=E(b)=j-k$. Considering the cases
$(j-k)$ even and $(j-k)$ odd, it is straight-forward to check that for
all $j,\, k\in\Bbb N$ with $j-k=0,1,\dots,n-2$, we have
\begin{equation}
N_2 =j-k.
\end{equation} 
Equations (3.11), (3.12), (3.13) and (3.14) complete the proof of (ii).
\end{itemize}
\medskip

By virtue of Theorem 3.3 and the identities (2.1), (2.2) and (3.8), it is
easy to see that we only have one possibility left that has not
been analyzed, namely,
\begin{equation}
1-n<c-b<0, \quad 1-n<b<0, \quad 1-n<c<0.
\end{equation} 
\bigskip

\noindent{\bf Theorem 3.4. } {\it Let $F=F(-n,b;\, c;\, z)$ where
$b$ and $c$ satisfy condition (3.15). If $-j<b<-j+1$, $j=1,\dots,n-1$;
$-k<c<-k+1$, $k=1,\dots,n-1$ and $-\ell<c-b<-\ell+1$, $\ell=1,\dots,n-1$,
then $F$ has no real zeros if $n+j+\ell$, $k+\ell$, $j+k$ are even,
one real zero in $(1,\infty)$ if $n+j+\ell$ is odd, one real zero in
$(0,1)$ if $k+\ell$ is odd and one real negative zero if $j+k$ is odd.}
\bigskip

\noindent{\bf Proof. } Under the restrictions (3.15), we have
\[
|1-c|=1-c, \quad |b+n|=b+n, \quad |b-c-n|=c-b+n.
\]
Then from (3.2), (3.3) and (3.4),
\[
X=E(1-c-n), \quad Y+E(b), \quad Z=E(c-b),
\]
and it follows from (3.15) that $X=Y=Z=0$. Also, $\sign
{-b\choose n}=(-1)^{n-j}$, $\sign {-c\choose n}=
(-1)^{n-k}$ and $\sign {b-c\choose n} =(-1)^{n-\ell}$.
The stated result then follows immediately from (3.5), (3.6) and (3.7).
\bigskip

\noindent{\bf Remark. }  We have not considered the asymptotic zero
distribution as $n \to \infty$ of $F(-n,b;\, c;\, z)$. There are recent interesting
results in this regard using different approaches, namely complex
analysis techniques [10], matrix theoretic tools [11], asymptotic analysis
of the Euler integral representation [3] and analysis of coefficients [8].

\bigskip

\begin{tabbing}
e-mail addresses: \= 036KAD@COSMOS.WITS.AC.ZA \\
\> 036JORD@COSMOS.WITS.AC.ZA
\end{tabbing}

\end{document}